\documentclass[11pt, a4paper]{amsart}

\usepackage{amsfonts,amsmath,amssymb, amscd,fullpage,mathtools}
\usepackage[all]{xy}
\usepackage{enumitem}
\usepackage{xcolor}

\newtheorem{theorem}{Theorem}[section]
\newtheorem{lemma}[theorem]{Lemma}
\newtheorem{definition}[theorem]{Definition}

\newtheorem{corollary}[theorem]{Corollary}

\theoremstyle{definition}
\newtheorem{rem}[theorem]{Remark}

\newtheorem{remark}[theorem]{Remark}
\newtheorem{example}[theorem]{Example}
\newtheorem*{claim}{Claim}

\newcommand\pf{\begin{proof}}
\newcommand\epf{\end{proof}}

\newcommand\eps{\varepsilon}

\newcommand\co{\operatorname{co}}

\DeclareMathOperator{\ad}{ad}

\DeclareMathOperator{\Ker}{Ker}

\DeclareMathOperator{\id}{id}

\numberwithin{equation}{section}

\newcommand\op{\operatorname{op}}
\newcommand\N{\mathbb{N}}
\newcommand\ot{\otimes}
\DeclareMathOperator{\Aut}{Aut}
\newcommand\cP{\mathcal{P}}
\newcommand\inv{^{-1}}
\DeclareMathOperator{\can}{can}

\def\lg{\langle}
\def\rg{\rangle}
\def\sl{\mathfrak{sl}}

\title{Hopf-Galois structures on ambiskew polynomial rings}

\author{Julien Bichon}
\address{Julien Bichon:	
Laboratoire de Math\'ematiques Blaise Pascal,
Universit\'e Clermont Auvergne \& CNRS,
Complexe universitaire des C\'ezeaux,
63178~Aubi\`ere Cedex, France}
\email{julien.bichon@uca.fr}

\author{Agust\'in Garc\'ia Iglesias}
\address{Agust\'in Garc\'ia Iglesias: Facultad de Matem\'atica, Astronom\'ia, F\'isica y Computaci\'on (FaMAF), 
	Medina Allende S/N, 
	Universidad Nacional de C\'ordoba,
	Ciudad Universitaria, C\'ordoba (X5000HUA),
	Argentina. }
\email{aigarcia@famaf.unc.edu.ar}

\begin{document}

\begin{abstract}
 We provide necessary and sufficient conditions to extend the Hopf-Galois algebra structure on an algebra $R$ to a generalized ambiskew ring based on $R$, in a way such that the added variables for the extension are skew-primitive in an appropriate sense. We show that the associated Hopf algebra is again a a generalized ambiskew ring, based on a suitable Hopf algebra $\underline{H}(R)$. Several examples are examined, including the Hopf-Galois objects over $U_q(\mathfrak{sl}_2)$.
\end{abstract}

\maketitle

\section{Introduction}\label{sec:intro}

	Let $R$ be an algebra over a field $k$ and consider the following data:
	\begin{enumerate}
		\item A pair of commuting  algebra automorphisms $\tau,\omega\in \Aut R$; set $\sigma\coloneqq \tau\omega=\omega\tau$.
		\item A $\sigma$-central element $c$; that is $c\in Z_\sigma(R)\coloneqq\{x\in R: xr=\sigma(r)x, \ r\in R\}$.
		\item A scalar $\xi \in k^*$.
	\end{enumerate}
	The associated \emph{generalized ambiskew polynomial algebra} $A=A(R,X,Y, \tau,\omega, c, \xi)$, see \cite{jor-well}, is the quotient of the free product $R*k\langle X,Y\rangle$ by the relations
	\begin{align*}
	Xr &= \tau(r)X, & Yr &= \omega(r)Y, & XY-\xi YX&=c, 
	\end{align*}
	for any $r \in R$. 
	
When $c$ is assumed to be central (so that $\sigma=\id$ or $\omega=\tau^{-1}$)	then $A=A(R,X,Y, \tau,\tau^{-1}, c, \xi)$ is called an \emph{ambiskew polynomial algebra} \cite{jor00}; we set  $A(R,X,Y, \tau, c, \xi)\coloneqq A(R,X,Y, \tau,\tau^{-1} c, \xi)$.

A common question is about the ring-theoretical properties such as simplicity, or the representation theory, of the resulting algebra $A$, see for example \cite{jor}, \cite{jor-well}, a useful summary of known results being given in \cite[Section 6]{brma}. 

Another interesting question is to find necessary and sufficient conditions on these data so that $A$ preserves some extra structure on $R$. 
When $R$ is a Hopf algebra, Brown-Macaulay \cite{brma} discuss the case in which the Hopf structure on $R$ extends to a Hopf structure on  the ambiskew polynomial Hopf algebra
$A_R=A(R,X, Y, \tau, c, \xi)$ via the formula
\begin{align}
	\Delta(X)&=X \otimes 1 + g \otimes X, & \Delta(Y)&=Y \otimes 1 + h \otimes Y
\end{align}
for some group-like elements $g,h \in G(R)$; here, $c$ is central\footnote{In loc.cit.~the authors choose the notation $(X_+,X_-)$ and $(g_+,g_-)$ instead of $(X,Y)$ and $(g,h)$, respectively; and use $h$ instead of $c$, $\sigma$ instead of $\tau$.}. This generalized previous results from \cite{har} in the case when  $R$ is commutative.
As an example, $U_q(\mathfrak{sl}_2)$ can be constructed from $k\mathbb Z$ by this process.

%
%
%\subsection{Ambiskew Polynomial rings}
%
%Let $R$ be an algebra and consider the following data:
%\begin{enumerate}
%\item An algebra automorphism $\sigma$ of $R$,
%\item an element $c\in Z(R)$, 
%\item a scalar $\xi \in k^*$.
%\end{enumerate}
%The associated ambiskew polynomial algebra $A=A(R,X_+,X_-, \sigma, c, \xi)$ is the quotient of the free product $R*k\langle X_+,X_-\rangle$ by the relations
%\begin{align*}
%X_+r &= \sigma(r)X_+, & X_-r &= \sigma^{-1}(r)X_-, & X_+X_--\xi X_-X_+&=c, 
%\end{align*}
%for any $r \in R$. 

In this article, we generalize these results in two directions:
\begin{enumerate}[label=(\Alph*)]
\item\label{itemA} We discuss the extension of the Hopf algebra structure on $R$ to a generalized ambiskew polynomial algebra $A$. 
\item\label{itemB} We further extend this to the case when $R$ is a Hopf-Galois object over some Hopf algebra. 
\end{enumerate}

Hopf-Galois objects are natural generalizations of Hopf algebras with a Galois-theoretic flavour, introduced in \cite{kt}, and have proved to be fundamental tools in analysing tensor categories of comodules \cite{ulb,sc1} as well as in recent classification questions for pointed Hopf algebras \cite{aagmv,ag}.

While our first goal was to extend the results of \cite{brma} to Hopf-Galois objects in the case of ambiskew polynomial rings, our analysis evidenced that the proper context for doing that was that of generalized ambiskew polynomial rings, leading to the results in item (A).

The punchline for  item (B) is that while having to deal simultaneously with $R$ and the Hopf algebra coacting on it  can look difficult, it turns out that, fortunately, results of Grunspan \cite{gru} and Schauenburg \cite{scqt} ensure that there is exists a definition of Hopf-Galois objects that does not use any Hopf algebra, in terms of an algebra map
\[
\mu\colon R\to R\ot R^{\op}\ot R,
\]
subject to certain axioms, see Section \ref{sec:hopf-galois}. We start with a Hopf-Galois algebra $R$, and discuss the ingredients with  necessary and sufficient conditions to extend the Hopf-Galois structure of $R$ to $A(R,X, Y, \sigma, c, \xi)$ via the formula 
\begin{align}
\begin{split}\label{eqn:muY}
\mu(X)&=X\otimes 1 \otimes 1 - g \otimes g^{-1}X\otimes 1 + g\otimes g^{-1}\otimes X,\\
\mu(Y)&=Y \otimes 1 \otimes 1 - h\otimes h^{-1}Y\otimes 1 + h\otimes h^{-1}\otimes Y,
\end{split}
\end{align}
for some quasi-central group-like elements $g,h \in G_0(R)$, see Definition \ref{def:group-like}. %If $t\in G_0(R)$, then we write $t\cdot r\coloneqq trt^{-1}$, $r\in R$.

Given a Hopf-Galois algebra $(R,\mu)$ as above, there exists  a canonical Hopf algebra $\underline{H}(R)$ such that $R$ is  an $\underline{H}(R)$-Galois object \cite{gru,scqt}. As for the first item \ref{itemA}, the need of a generalization of the results in \cite{brma} becomes evident when one considers the Hopf algebra associated to the Hopf-Galois algebra $A_R=A(R,X, Y, \tau, c, \xi)$ extending a Hopf-Galois algebra $R$: this Hopf algebra is an extension of $\underline{H}(R)$ not in terms of an ambiskew polynomial algebra but of a generalized ambiskew polynomial algebra $A$. See  Corollary \ref{coro:notid} and Subsection  \ref{sec:illustration} for an example.

The paper is organized as follows. In Section \ref{sec:hopf-galois} we recall the definition of Hopf-Galois algebras $(R,\mu)$ and their realization as Hopf-Galois objects over a Hopf algebra $\underline{H}(R)$; we discuss the concepts of group-like and skew-primitive elements in this setting. Section \ref{sec:ambiskew-hopf} is devoted to the answer of item \ref{itemA} above: in Theorem \ref{thm:hopf-algebra} we find necessary and sufficient conditions so that a generalized ambiskew polynomial algebra $A_R$ extends the structure of a Hopf algebra $R$ in a way such that the added variables are skew-primitive; we recover some of the results in \cite{brma} as a corollary. Next, in Section  \ref{sec:ambiskew-hopfgalois} we find necessary and sufficient conditions so that a generalized ambiskew polynomial algebra extends the structure of a Hopf-Galois algebra $R$ as in (\ref{eqn:muY}) and we provide an answer for item \ref{itemB} in Theorem \ref{thm:mainmain}. In Section \ref{sec:hopf-to-ambiskew} we show in Theorem \ref{thm:H(A)} that the Hopf algebra associated to a Hopf-Galois ambiskew polynomial algebra $A_R$ is a Hopf ambiskew  polynomial algebra, associated to $\underline{H}(R)$. We conclude the article with a series of examples in Section \ref{sec:examples}.

\subsection*{Notation and conventions}
We work over a field $k$. All algebras, vector spaces and unadorned tensor products $\ot$ are assumed to be defined over $k$. We assume that the reader has familiarity with Hopf algebra theory, for which the book \cite{mon} is a convenient reference. We use standard notations, in particular $\Delta$, $\varepsilon$ and $S$ stand respectively  for the comultiplication, counit and antipode of a Hopf algebra, and 
Sweedler's notation $h\mapsto h_{(1)}\ot h_{(2)}$, resp.~$m\mapsto m_{(0)}\ot m_{(1)}$ for the comultiplication, resp.~the coaction on a (right) comodule $M$, of a Hopf algebra $H$; $h\in H, m\in M$

\subsection*{Acknowledgements}
The authors were supported by CONICET, FONCyT PICT-2015-2845, Secyt (UNC), the MathAmSud project GR2HOPF, and the Emergence Program of the I-SITE project CAP 20-25 of Clermont Auvergne University.

\section{Hopf-Galois algebras}\label{sec:hopf-galois}

\subsection{Basic definitions}
The definition below was introduced by Grunspan \cite{gru} with more axioms, and simplified later by Schauenburg \cite{scqt}, under the name {\it quantum torsor}. In view of Theorem \ref{thm:hr} given in the next subsection, it seems natural to call such a structure a Hopf-Galois algebra. 

\begin{definition}
 A Hopf-Galois algebra is a non-zero algebra $R$ together with an algebra map 
$$\mu : R \longrightarrow R \otimes R^{\rm op} \otimes R$$
such that 
\begin{align}
\label{eqn:galois1}(\mu\otimes {\rm id}_R\otimes {\rm id}_R)\circ \mu&= ({\rm id}_R\otimes {\rm id}_R\otimes \mu )\circ \mu,\\
\label{eqn:galois2}(m\otimes {\rm id}_R)\circ \mu&=\eta \otimes {\rm id}_R, \quad ({\rm id}_R\otimes m) \circ \mu= {\rm id}_R\otimes \eta, 
\end{align}
where $m: R \otimes R \to R$ and $\eta : k \to R$ denote the respective multiplication and unit of $R$.
\end{definition}
We shall write 
\[
\mu(r)=r_{(1)}\ot r_{(2)}\ot r_{(3)}, \quad r\in R.
\]
Hence, for any $r\in R$, 
\begin{align*}
r\ot 1&=r_{(1)}\ot r_{(2)}r_{(3)}, & 1\ot r&=r_{(1)}r_{(2)}\ot r_{(3)};
\end{align*}
and we may set
\[
r_{(1)}\ot r_{(2)}\ot r_{(3)}\ot r_{(4)}\ot r_{(5)}\coloneqq\mu(r_{(1)})\ot r_{(2)}\ot r_{(3)}=r_{(1)}\ot r_{(2)}\ot \mu(r_{(3)}).
\]

\begin{example}
The basic example is a Hopf algebra $H$  with 
\[\mu(x)=x_{(1)}\otimes S(x_{(2)}) \otimes x_{(3)},  \ x\in R.\] We denote by $\underline{R}(H)$ the resulting Hopf-Galois algebra. This defines a functor from the category of Hopf algebras to the category of Hopf-Galois algebras (a morphism of Hopf-Galois algebras being an algebra map commuting with $\mu$ in the obvious way).

A Hopf-Galois algebra arises from a Hopf algebra as above if and only if there exists an algebra map $\alpha : R \to k$. Indeed, starting with $R$ and such an $\alpha$, the Hopf algebra $H$ whose counit is $\alpha$ and with the other  structure maps  defined by
$$\Delta(x)=\alpha(x_{(2)}) x_{(1)}\otimes x_{(3)}, \  S(x)=\alpha(x_{(1)}x_{(3)})x_{(2)}$$
is such that $R=\underline{R}(H)$.
\end{example}

\begin{example}\label{ex:Weyl}
 A more significant example is the Weyl algebra $$A_1(k)=k\langle x,y \ | \ xy-yx=1\rangle$$ with Hopf-Galois structure defined by
\begin{align*}
\mu(x)&= x \otimes 1 \otimes 1 - 1 \otimes x \otimes 1 + 1 \otimes 1 \otimes x, &  \mu(y)&= y \otimes 1 \otimes 1- 1 \otimes y \otimes 1 + 1 \otimes 1 \otimes y.
\end{align*}
\end{example}

More examples will be discussed in the final section. We now introduce group-like elements and skew-primitives in Hopf-Galois algebras.
%The above definition was introduced by Grunspan \cite{gru} with more axioms, and simplified later by Schauenburg \cite{scqt}, under the name {\it quantum torsor}. In view of their result below, it seems natural to call this structure a Hopf-Galois algebra.

\begin{definition}\label{def:group-like}
A group-like element in a Hopf-Galois algebra  $R$ is an invertible element $g\in R$ such that
\begin{align}\label{eqn:grouplike}
\mu(g)=g\ot g^{-1}\ot g.
\end{align}
We denote by $G_0(R)$ the set of group-like elements in $R$; $g\in G_0(R)$ is quasi-central if there exists a character $\alpha$ on $G_0(R)$ such that 
\[
xg = \alpha(x)gx,  \text{ for any } x\in G_0(R).
\]
\end{definition}
If $g\in G_0(R)$, then we write $g\cdot r\coloneqq grg^{-1}$.

\begin{remark}\label{rem:gp}
It should be noticed that, contrary to the Hopf algebra case, if $g\in R$ is group-like, then $\lambda\,g$ is grouplike for any nonzero scalar $\lambda\in k$. If we consider $G_0(R)$ acting on $R$ by conjugation, then $g$ and $\lambda\,g$ define the same operator.
\end{remark}

\begin{definition}
Let $R$ be a Hopf-Galois algebra and let $g,h\in G_0(R)$. The set $\cP_{g,h}^0(R)$ of $(g,h)$-skew primitive elements is the collection of those 
$x\in R$ such that 
\[
\mu(x)=x\ot h^{-1}\ot h-g\ot g^{-1}xh^{-1}\ot h+g\ot g^{-1}\ot x.
\]
\end{definition}

\begin{rem}\label{rem:prim} If $x\in \cP^0_{g,h}(R)$, then $h^{-1}x\in \cP^0_{h^{-1}g,1}(R)$ and $xg^{-1} \in \cP_{1,hg^{-1}}^0(R)$.
Notice also that $g,h\in \cP_{g,h}^0(R)$. Indeed,
\[
\mu(h)=h\ot h^{-1}\ot h=h\ot h^{-1}\ot h-g\ot g^{-1}hh^{-1}\ot h+g\ot g^{-1}\ot h.
\]
and similarly for $g$.
\end{rem}

%\subsubsection{On extra axioms: the Grunspan map $\theta$} Do we us this ?

%Consider the map $\theta:R\to R$ given by $\theta(r)=r_{(1)}r_{(2)(3)}r_{(2)(2)}r_{(2)(1)}r_{(3)}$, $r\in R$. Then
%\begin{align}
%\label{eqn:gru1}\mu\theta&=\theta^{\ot3}\mu, \quad \text{and}\\
%\label{eqn:gru2}(\id\ot\mu^{\op}\ot \id)\mu&=(\id\ot\id\ot\theta\ot\id\ot\id)(\mu\ot\id\ot\id)\mu.
%\end{align}
%The map $\theta$ is called a Grunspan map and, as we see, it is univocally determined by $\mu$. The existence of a map $\theta$ satisfying \eqref{eqn:gru1} and \eqref{eqn:gru2} was originally stated as an axiom in \cite{gru}.

%In Sweedler's type notation, \eqref{eqn:gru2} reads:
%\begin{align*}
%r_{(1)}\ot r_{(2)(3)}\ot r_{(2)(2)}\ot r_{(2)(1)}\ot r_{3}=r_{(1)}\ot r_{(2)}\ot \theta(r_{(3)})\ot r_{(4)}\ot r_{(5)}.
%\end{align*}
%In particular, we shall make use of the following derivation of \eqref{eqn:gru2}:
%\begin{multline}\label{eqn:gru3}
%r_{(1)}\ot r_{(2)}\ot r_{(3)}\ot r_{(4)(3)}\ot r_{(4)(2)}\ot r_{(4)(1)}\ot r_{(5)}\\
%=r_{(1)}\ot r_{(2)}\ot r_{(3)}\ot r_{(4)}\ot \theta(r_{(5)})\ot r_{(6)}\ot r_{(7)}.
%\end{multline}

\subsection{From Hopf-Galois algebras to Hopf algebras}\label{sec:HR}
We recall that a (right) Hopf-Galois object $R$ over a Hopf algebra $H$ is a (right) comodule algebra $R$ over $H$ with trivial coinvariants $R^{\co H}=k$ and such that the 
{\it Galois map} 
\[
\can\colon R\ot R\to R\ot H, \qquad \can(r\ot s)=rs_{(0)}\ot s_{(1)}, \ r,s\in R
\]
is bijective.

The following result relates what we called Hopf-Galois algebras and classical Hopf-Galois objects.

\begin{theorem}[Grunspan-Schauenburg, \cite{gru,scqt}]\label{thm:hr}
Let $R$ be an algebra. The following assertions are equivalent.
 \begin{enumerate}
  \item There exists an algebra map $\mu : R \longrightarrow R \otimes R^{\rm op} \otimes R$ making $R$ into a Hopf-Galois algebra.
\item $R$ is a right Hopf-Galois object over some Hopf algebra $\underline{H}(R)$.
 \end{enumerate}
\end{theorem}
Therefore the notion of Hopf-Galois algebra provides a definition of Hopf-Galois object which is Hopf algebra free. %We will provide more details on this construction on Section \ref{sec:HR}.

The constuction of $\underline{H}(R)$ in the theorem is as follows. Let $R=(R,\mu)$ be a Hopf-Galois algebra. Then $R$ is a right Hopf-Galois object for the Hopf algebra $\underline{H}(R)$ defined as the subalgebra
\begin{equation}\label{eqn:HR}
\underline{H}(R)=\left\{ x \ot y  :  x  y_{(1)}\ot y_{(2)}\ot y_{(3)}= 1\ot x \ot y \right\}\subseteq R^{\op}\ot R
\end{equation}
with coalgebra structure
\begin{equation}\label{eqn:HR-coalg}
\Delta( x \ot y  )=  x  \ot y_{(1)}\ot y_{(2)}\ot y_{(3)}, \quad \eps( x \ot y )= x y .
\end{equation}
The coaction is given by $\mu$, as $\mu\colon R\to R\ot \underline{H}(R)$.
\begin{rem}\label{rem:sweedler}
	Here, we follow Schauenburg in \cite[\S 3]{scqt} and we write a generic element $\sum x_i\ot y_i\in R^{\op}\ot R$ as a single tensor $x\ot y$, ``in the spirit of Sweedler's notation''. Definition \eqref{eqn:HR} and formula \eqref{eqn:HR-coalg} should be interpreted in this sense.
\end{rem}
Conversely, let $H$ be a Hopf algebra and let $R$ be a right Hopf-Galois object, with coaction $r\mapsto \rho(r)=r_{(0)}\ot r_{(1)}$, $r\in R$. Then $R$ is a Hopf-Galois algebra with
\[
\mu(r)=r_{(0)}\ot \can^{-1}(1\ot r_{(1)}).
\]

\medskip

We now use the construction $\underline{H}(R)$ to relate the group-likes and skew-primitives of the previous subsection with ordinary group-likes and skew-primitives in a Hopf algebra.

\begin{lemma}\label{lem:group}
Let $R$ be a Hopf-Galois algebra.
	\begin{enumerate}
		\item There is a short exact sequence of groups
		\begin{align*}
		1\to  k^\ast \rightarrow G_0(R)\stackrel{\varphi}{\longrightarrow} G(\underline{H}(R))\to 1,
		\end{align*}
		where $\varphi(g)=g^{-1}\ot g$.
		\item For any $g\in G_0(R)$, there is a linear short exact sequence
		\begin{align*}
		0\to  k \rightarrow \cP_{g,1}^0(R) \stackrel{\psi}{\longrightarrow} \cP_{\varphi(g),1}(\underline{H}(R))\to 0, 
		\end{align*}
		where $\psi(x)=g^{-1}\ot x-g^{-1}x\ot 1$.
		%The inverse is given by 
		%\begin{align}\label{eqn:map-primitive-inverse}
		%x\ot y\mapsto (y-\beta)((1-g)x+\alpha\,g)
		%\end{align}
		%for certain scalars $\alpha,\beta\in k$, see \eqref{eqn:scalars}.
	\end{enumerate}
\end{lemma}
\pf
(1) On the one hand, we see that $k^\ast\subset G_0(R)$, as in Remark \ref{rem:gp}. 
If $g\in G_0(R)$, then it is clear that $g^{-1}\ot g\in G(\underline{H}(R))$.
Moreover, $\varphi(\lambda)=1\ot 1$ for any $\lambda\in k^\ast$.

Now, if $u=x\ot y=\sum_i x_i\ot y_i\in G(\underline{H}(R))$, then $\eps(x\ot y)=xy=\sum_ix_iy_i=1$ and
\[
\Delta(x\ot y)=x\ot y\ot x\ot y=\sum_{i,j} x_i\ot y_i\ot x_j\ot y_j.
\]
On the other hand, by definition 
\[
\Delta(x\ot y)=x\ot y_{(1)}\ot y_{(2)}\ot y_{(3)}.
\]
Hence, by multiplying the factors in the middle, we get
\[
x\ot 1\ot y=\sum_i x_i\ot 1\ot y_i=\sum_{i,j} x_i\ot y_ix_j\ot y_j.
\]
If we assume that $\{x_i\}_i$ is linearly independent, we obtain, for any $i$, that 
$1\otimes y_i=\sum_jy_ix_j\otimes y_j$ and hence \[1\otimes \mu(y_i)=\sum_jy_ix_j\otimes \mu(y_j).\]
Multiplying the left handed tensors, this gives 
\[\mu(y_i)=\sum_jy_ix_jy_{j(1)}\otimes y_{j(2)}\otimes y_{j(3)}=\sum_jy_i\otimes x_j\otimes y_j=y_i\otimes u.\]
Applying the same reasoning to $u^{-1}=\sum_kx'_k\otimes y'_k$ yields
\[\mu(y'_k)=y'_k\otimes u^{-1}\]
for any $k$. Hence for any $i,k$, we have $\mu(y_iy'_k)=y_iy'_k\otimes 1\otimes 1$ and
$\mu(y'_ky_i)=y'_ky_i\otimes 1\otimes 1$, so that there exists scalars $\lambda_{i,k}$, $\mu_{k,i}$ such that
 $y_iy'_k=\lambda_{i,k}$ and $y_k'y_i=\mu_{k,i}$. There exists $j,l$ such that $y_jy'_l\not =0$, so we see from the previous identities that $\lambda_{j,l}=\mu_{l,j}$ and that $y'_l$ is invertible in $R$. Hence for any $i$ we have $y_i=\lambda_{i,l}(y_l')^{-1}$. We get $u=(\sum_i \lambda_{i,l}x_i)\otimes (y'_l)^{-1}=g^{-1}\otimes g$ for $g={y'_l}^{-1}$, and we then see that $g \in G_0(R)$.
%\in R^{{\rm co}\underline{H}(R)}$ with $R^{{\rm co}\underline{H}(R)}=k1$ since $R$ is a Galois object over $\underline{H}(R)$.
%$y_ix_j=\delta_{i,j}$ and thus $x\ot y=\sum_i y_i^{-1}\ot y_i$. As $xy=1$, we see that $x\ot y=g^{-1}\ot g=\varphi(g)$, for some $g\in G_0(R)$.

\medbreak

(2) We know that $k\subset \cP_{g,1}^0(R)$, see Remark \ref{rem:prim}. Also, it is easy to check that $g^{-1}\ot x-g^{-1}x\ot 1\in P_{\varphi(g),1}(\underline{H}(R))$ if $x\in \cP_{1,g}^0(R)$, and that $\Ker(\psi)=k$.

Now, let $u=x\ot y=\sum x_i\ot y_i\in \cP_{\varphi(g),1}(\underline{H}(R))$. We assume that $\{y_i\}$ is linearly independent.
We have 
$xy_{(1)}\ot y_{(2)}\ot y_{(3)}=1\ot x\ot y$ and
\begin{equation}\label{expression}
x\ot y_{(1)}\ot y_{(2)}\ot y_{(3)}= x\ot y\ot 1\ot 1+g^{-1}\ot g\ot x\ot y
\end{equation}
so that multiplying the factors in the middle as before:
\[
x\ot 1\ot y= x\ot y\ot 1+g^{-1}\ot gx\ot y, \ i.e. \] 
\[\sum_i x_i\ot 1\ot y_i= \sum_i x_i\ot y_i\ot 1+\sum_ig^{-1}\ot gx_i\ot y_i. %\qquad (=g^{-1}\ot 1\ot h-g^{-1}h\ot1\ot 1). 
\]
If $1\not \in {\rm Span}\{y_i\}$, we get that $u=\sum_ix_i\ot y_i=0$. Otherwise, there exists for any $i$, a scalar $\lambda_i$ such that
\[ x_i\otimes 1 = g^{-1}\otimes gx_i + \lambda_i u
\]
with $\lambda_k\not=0$ for some $k$. This gives \[u=g^{-1}\otimes z-g^{-1}z\otimes 1\] for $z=-\lambda_k^{-1}gx_k$. %If $z=\alpha g$ for some scalar $\alpha$, then  $u= \alpha (g^{-1}\ot g -1\otimes 1)=\psi(\alpha(g-1))$ with $\alpha(g-1)\in  \cP_{g,1}^0(R)$ by Remark \ref{rem:prim}. We then assume that $z$ and $g$ are linearly independent. 
Going back to the expression \eqref{expression}, we get
\[g^{-1}\otimes \mu(z) -g^{-1}z\otimes 1 \ot 1 \ot 1 = g^{-1}\ot z \otimes  1\ot 1 -g^{-1}z\otimes 1 \ot 1\ot 1+ g^{-1}\otimes g \ot g^{-1}\ot z - g^{-1}\otimes g\ot g^{-1}z\ot 1.\]
This shows that $z\in  \cP_{g,1}^0(R)$, and we have $u=\psi(z)$. \epf

\section{Ambiskew Hopf algebras}\label{sec:ambiskew-hopf}

In this section we assume that $R$ is a Hopf algebra and consider the generalized ambiskew algebra $A=A(R,E,F,\tau,\omega, c, \xi)$ as in the introduction. We look for necessary and sufficient conditions to extend the comultiplication in $R$ to $A$ via the formulas
\begin{align} \label{eqn:DeltaE}
\Delta(E)&=E \otimes 1 + g \otimes E, & \Delta(F)&=F \otimes 1 + h \otimes F,
\end{align}
for some (necessarily group-like) elements $g,h\in G(R)$.
 
When the element $c$ is assumed to be central, or equivalently when $\omega=\tau^{-1}$ (so $\sigma=\id$), this question has been addressed in \cite{brma}. Our result reads as follows.
 
\begin{theorem}\label{thm:hopf-algebra}
Let $H$ be a Hopf algebra, let $\tau,\omega\in\Aut H$ be such that $\tau\omega=\omega\tau$  and define $\sigma=\tau\omega$. Let $c\in Z_\sigma(H)$ be $\sigma$-central and $\xi \in k^*$. 
Then the generalized ambiskew polynomial algebra $A(H,E,F, \tau,\omega, c, \xi)$ has a Hopf algebra structure extending that of $H$ and such that \eqref{eqn:DeltaE} holds for some $g,h\in G(H)$ if and only if 
\begin{itemize}
\item there are characters $\alpha,\beta$ of $H$, with $\alpha\ast\beta=\beta\ast\alpha$ and 
\begin{align}\label{eqn:tau-omega-hopf}
\tau(r)&=\alpha(r_{(1)})r_{(2)}, & \omega(r)=\beta(r_{(1)})r_{(2)}.
\end{align}
\begin{align}\label{eqn:alpha-beta}
\alpha(r_{(1)})r_{(2)}&=g\cdot r_{(1)}\alpha(r_{(2)}), 
&
\beta(r_{(1)})r_{(2)}&=h\cdot r_{(1)}\beta(r_{(2)})
\qquad r\in H,
\end{align}
\item the group-like elements $g,h$ are central in $G(H)$ and
$\alpha(h)=\beta(g)^{-1}=\xi$,
\item  the $\sigma$-central element $c$ is $(gh,1)$-skew primitive.
\end{itemize}
\end{theorem}
\pf
Assume $A=A(H,E,F, \tau,\omega, c, \xi)$ is a Hopf algebra as in the statement. Then, by comparing $\Delta(Xr)$ and $\Delta(\tau(r)X)$ we obtain that
\[
\tau(r_{(1)})\ot r_{(2)}=\tau(r)_{(1)}\ot \tau(r)_{(2)}=g\cdot r_{(1)}\ot \tau(r_{(2)}).
\]
Hence, if $\alpha=\eps\circ\tau:H\to k$ we get that $\alpha(r_{(1)})r_{(2)}=g\cdot r_{(1)}\alpha(r_{(2)})$ and $\tau(r)=\alpha(r_{(1)})r_{(2)}$. Similarly, 
by comparing $\Delta(Yr)$ and $\Delta(\omega(r)Y)$ we obtain that $\omega(r)=\beta(r_{(1)})r_{(2)}$ for $\beta=\eps\circ\omega$ and both identities in \eqref{eqn:alpha-beta} hold. Also, since $\tau\omega=\omega\tau$, it follows that $\alpha\ast\beta=\beta\ast\omega$.

Now, observe that \eqref{eqn:alpha-beta} applied to $r=h$ gives $g\cdot h=h$, that is $gh=hg$. We use this to see that
\begin{align*}
\Delta(EF-\xi FE)=&(EF-\xi FE)\ot 1 + gh\ot (EF-\xi FE)\\
&\qquad +(\alpha(h)-\xi)hE\ot F+(1-\xi\beta(g))gF\ot E\\
&=\Delta(c)\in H\ot H
\end{align*}
and thus, since $A$ is a free $R$-module with basis $\{X^aY^b:a,b\geq 0\}$, we have $\alpha(h)=\beta(g)^{-1}=\xi$ and $\Delta(c)=c\ot 1 + gh\ot c$, so $c=EF-\xi FE$ is a $\sigma$-central $(gh,1)$-skew primitive element. 
The previous computations allow as well to prove converse statement.
\epf

\begin{remark}\label{rem:either}
	Notice that $g\cdot c=g\sigma(g)^{-1}c=\alpha(g)^{-1}\beta(g)^{-1}c=\xi\alpha(g)^{-1}c$. We can thus compute $\tau(c)$ in two ways and get
	\begin{align*}
	\xi(\alpha(g)^{-1}-\alpha(g))c=\alpha(c)(1-gh).
	\end{align*}
	So either $\alpha(g)=\pm1$ (and $\alpha(c)=0$ or $1=gh$) or $c=\frac{\alpha(ghc)}{1-\alpha(g)^2}(1-gh)$.
	
	Similarly, if we compute $\omega(c)$, we get that
	\begin{align*}
	\xi^{-1}(\beta(h)^{-1}-\beta(h))c=\beta(c)(1-gh),
	\end{align*}
	that is $\beta(h)=\pm1$ (and $\beta(c)=0$ or $1=gh$) or $c=\frac{\beta(ghc)}{1-\beta(h)^2}(1-gh)$.
	\end{remark}

Theorem \ref{thm:hopf-algebra}  extends the result in \cite{brma} as mentioned; we recall this here as a corollary. 

\begin{corollary}\cite{brma}\label{cor:hopf-algebra}
Assume $c$ is central (i.e.~$\omega=\tau^{-1}$ so $\sigma=\id$). Then the ambiskew polynomial algebra $A(R,E,F, \tau, c, \xi)=A(R,E,F, \tau,\tau^{-1}, c, \xi)$ has a Hopf algebra structure extending that of $H$ and such that \eqref{eqn:DeltaE} holds for some $g,h\in G(H)$ if and only if there are 
\begin{itemize}
\item a character $\chi$ of $R$,
\item central group-like elements $g,h\in Z(G(R))$ such that $gh\in Z(R)$ and
\[\chi(g)=\chi(h)=\xi,\]
\item a central $(gh,1)$-skew primitive element $c$,
\end{itemize}
such that the following\footnote{In the notation from \cite{brma}, this condition reads $\tau^l_\chi=\ad_l(y_+)\circ \tau^r_\chi$ and then $\tau=\tau^l_\chi$.} holds:
\begin{align*}
\chi(r_{(1)})r_{(2)}=g\cdot r_{(1)}\chi(r_{(2)}), \qquad r\in R.
\end{align*}
In this setting,
\begin{align*}
\tau(r)=\chi(r_{(1)})r_{(2)}=g\cdot r_{(1)}\chi(r_{(2)}).
\end{align*}
Moreover, 
\begin{align}\label{eqn:br-either}
\begin{split}
\text{either } & \bullet\ \xi=\pm1; \ (\text{and }gh=1 \text{ or } \chi(c)=0),\\
\text{or else } & \bullet\ \xi\neq \pm1 \ \text{and } c=\lambda(1-gh), \text{ for }\lambda=\frac{\chi(c)}{1-\xi^2}\in k.\end{split}
\end{align}
\qed
\end{corollary}
\begin{rem}
Conditions \eqref{eqn:br-either} are stated in \cite[\S3.1]{brma}, in this setting they follow by Remark \ref{rem:either}.
\end{rem}

%
%
%\subsection{A Yetter-Drinfeld approach}
%
%Let $R$ be a Hopf algebra. Recall that the category $\ydr$ of (left) Yetter-Drinfeld modules is the category with objects given by simultaneously $R$-modules and $R$-comodules $M$ such that the action $\cdot\colon R\ot M\to M$ and the coaction $\lambda\colon M\to R\ot M$ are compatible in the following sense:
%\[
%\delta(r\cdot m)=r_{(1)}m_{(-1)}\sS(r_{(3)})\ot r_{(2)}\cdot m_{0}, \qquad r\in R,  m\in M.
%\]
%
%Let us set $T_\xi\coloneqq k_{\xi}[X,Y]$, that is the quantum plane generated by $X,Y$ subject to the $\xi$-commutation rule:
%\[
%XY=\xi YX.
%\]
%
%\begin{lemma}
%Assume that $A_R$ is a Hopf algebra, with $\tau,\omega$ as in \eqref{eqn:tau-omega-hopf}. 
%Then  $T_{\xi}$ is a Hopf algebra in $\ydr$ with $X,Y\in\cP(T_\xi)$.
%In particular, 
%\[
%T_{\xi}\# R\simeq A(R,X_{\pm},\tau,\omega,0,\xi)
%\]
%is a Hopf algebra. 
%%Moreover, when $\xi\neq\pm 1$,  $gr A(R,X_{\pm},\sigma,c,\xi)=T_{\xi}\# gr R$.
%\end{lemma}
%\pf
%Define an $\rightharpoonup\colon R\ot T_\xi\to T_\xi$ and a coaction $\lambda\colon T_\xi\to R\ot T_\xi$ via
%\begin{align*}
%r\rightharpoonup E&=\tau(r)\,E, & r\rightharpoonup F&=\omega(r)\,F, & E&\stackrel{\lambda}{\longmapsto}g\ot X,
%& F&\stackrel{\lambda}{\longmapsto}h\ot F,
%\end{align*}
%then the lemma follows.
%\epf

\section{Ambiskew Hopf-Galois algebras}\label{sec:ambiskew-hopfgalois}

Let $(R,\mu)$ be a Hopf-Galois algebra and let us fix $A=A(R,X,Y, \tau,\omega, c, \xi)$ an associated generalized ambiskew polynomial algebra. We investigate when $A$ is again a Hopf-Galois algebra, extending $R$.

\begin{theorem}\label{thm:mainmain}
	Let $R=(R,\mu)$ be a Hopf-Galois algebra. Assume given a sextuple $(\tau,\omega,g,h,c,\xi)$ where
	\begin{enumerate}
		\item $\tau,\omega\in\Aut(R)$ are commuting algebra automorphisms satisfying, for any $r\in R$,
		\begin{align}
		\label{eqn:id1} \tau(r)_{(1)}\ot \tau(r)_{(2)}\ot \tau(r)_{(3)}&=\tau(r_{(1)})\ot r_{(2)} \otimes r_{(3)}=\tau(r_{(1)})\ot \tau(r_{(2)}) \otimes \tau(r_{(3)} ),\\
		\label{eqn:id1w} \omega(r)_{(1)}\ot \omega(r)_{(2)}\ot \omega(r)_{(3)}&=\omega(r_{(1)})\ot r_{(2)} \otimes r_{(3)}=\omega(r_{(1)})\ot \omega(r_{(2)}) \otimes \omega(r_{(3)} ).
		\end{align}
		\item $g,h\in G_0(R)$ are quasi-central group-like elements such that, for any $r\in R$,
		\begin{align}
		\label{eqn:id2} g\cdot r_{(1)}\otimes g\cdot r_{(2)}\ot r_{(3)}&=\tau(r_{(1)})\ot \tau(r_{(2)}) \otimes r_{(3)},\\
		\label{eqn:id3} h\cdot r_{(1)}\otimes h\cdot r_{(2)}\ot r_{(3)}&=\omega(r_{(1)})\ot \omega(r_{(2)}) \otimes r_{(3)}. 
		\end{align}
		\item $\xi\in k^\ast$ is a scalar such that  $\tau(h)=\xi\, h$, $\omega(g)=\xi\inv\, g$.
		\item $c\in Z_{\sigma}(R)\cap \cP_{gh,1}(R)$, for $\sigma\coloneqq\tau\omega$; that is $c$ is a $\sigma$-central skew-primitive element:
		\[
		\mu(c)=c\ot1\ot1-gh\ot (gh)^{-1}c\ot 1+gh\ot (gh)^{-1}\ot c
		\]
		and $cr=\sigma(r)rc$, for all $r\in R$.
	\end{enumerate}
	Then there is a unique Hopf-Galois algebra structure on $A=A(R,X,Y,\tau,\omega,c,\xi)$ extending $\mu$ and such that 
\begin{align}
\begin{split}\label{eqn:muY-thm}
\mu(X)&=X\otimes 1 \otimes 1 - g \otimes g^{-1}X\otimes 1 + g\otimes g^{-1}\otimes X,\\
\mu(Y)&=Y \otimes 1 \otimes 1 - h\otimes h^{-1}Y\otimes 1 + h\otimes h^{-1}\otimes Y.
\end{split}
\end{align}
	Conversely, given commmuting $\tau,\omega\in\Aut(R)$, $c\in Z_\sigma(R)$ for $\sigma=\tau\omega$ and $\xi\in k^\ast$, if \eqref{eqn:muY-thm} defines a Hopf-Galois algebra structure on $A(R,X,Y,\tau,\omega,c,\xi)$ extending $\mu$, for some $g,h\in G_0(R)$, then conditions $(1)$ to $(4)$ are satisfied.
\end{theorem}

%\subsection{The proof}
\pf
We assume there is an extension $\mu:A\to A\ot A^{\op} \ot A$ and elements $g,h\in G(R)$ such that \eqref{eqn:muY-thm} holds. Recall the notation $x\cdot r\coloneqq xrx^{-1}$, for $x\in G(R)$, $r\in R$.

\smallskip

$\diamondsuit${\it The commutation rule.}
On the one hand, we have, in $A\ot A^{\op}\ot A$,
\[
\mu(Xr)=\mu(\tau(r)X)=(\tau(r)_{(1)}\ot \tau(r)_{(2)}\ot \tau(r)_{(3)})\mu(X).
\]
On the other,
\begin{align*}
	\mu&(Xr)=Xr_{(1)} \otimes r_{(2)} \otimes r_{(3)} - gr_{(1)} \otimes r_{(2)}g^{-1} X \otimes r_{(3)} + gr_{(1)} \otimes r_{(2)}g^{-1} \otimes Xr_{(3)}\\
	&=\tau(r_{(1)})X \otimes r_{(2)} \otimes r_{(3)} - (g\cdot r_{(1)}) g\otimes g^{-1} X \tau^{-1}(g\cdot r_{(2)})\otimes r_{(3)} 
	\\&\qquad + (g\cdot r_{(1)}) g\otimes g^{-1}(g\cdot r_{(2)}) \otimes \tau(r_{(3)})X\\
	&=(\tau(r_{(1)})\ot r_{(2)} \otimes r_{(3)} )(X \otimes 1 \otimes 1) - (g\cdot r_{(1)}\otimes \tau^{-1}(g\cdot r_{(2)})\ot r_{(3)}) (g\ot g^{-1} X \otimes 1) \\
	&\qquad + (g\cdot r_{(1)} \otimes  g\cdot r_{(2)}\otimes \tau(r_{(3)}))     (g\otimes g^{-1} \otimes X).
\end{align*}
Since $A$ is a free $R$-module with basis $\{X^aY^b:a,b\in\N\}$, we get the following identities:
\begin{align*}
	\tau(r)_{(1)}\ot \tau(r)_{(2)}\ot \tau(r)_{(3)}&=\tau(r_{(1)})\ot r_{(2)} \otimes r_{(3)}=g\cdot r_{(1)} \otimes  g\cdot r_{(2)}\otimes \tau(r_{(3)})\\
	g\cdot r_{(1)}\otimes g\cdot r_{(2)}\ot r_{(3)}&=\tau(r_{(1)})\ot \tau(r_{(2)}) \otimes r_{(3)} 
\end{align*}
We may combine these identities to show that $\mu\circ \tau=\tau^{\ot3}\circ \mu$:
\begin{align*}
	\tau(r_{(1)})\ot r_{(2)} \otimes r_{(3)}&=g\cdot r_{(1)} \otimes  g\cdot r_{(2)}\otimes \tau(r_{(3)})=\tau(r_{(1)})\ot \tau(r_{(2)}) \otimes \tau(r_{(3)} )
\end{align*}
and thus we see that our identities are then equivalent to \eqref{eqn:id1} and \eqref{eqn:id2}. 
Similarly, from the identity $Yr=\omega(r)Y$ we obtain \eqref{eqn:id1w} and   \eqref{eqn:id3}.

\smallskip

$\diamondsuit${\it The group characters.} We make the following
\begin{claim}
	There are characters $\alpha,\beta,\gamma,\delta\in\widehat{G_0(R)}$ such that for any $x \in G_0(R)$,
	\begin{align}
	\label{eqn:chars1} \tau(x)&=\alpha(x) x &  \omega(x)&=\beta(x) x,\\
	\label{eqn:chars2} \tau(x)&=\gamma(x)g\cdot x, &  \omega(x)&=\delta(x)h\cdot x.
	\end{align}
	Hence, $\alpha(h)\beta(g)=\gamma(h)\delta(g)$ and both $g,h$ are quasi-central group-like elements in $G_0(R)$ with
	\begin{align}\label{eqn:gh}
gh(hg)^{-1}&\in k^\ast \subset Z(R), &	(gh)\ot (gh)^{-1}&=(hg)\ot (hg)^{-1}.
	\end{align}
\end{claim}

As for the proof, notice that combining \eqref{eqn:galois2} and \eqref{eqn:id1} we get
\begin{align*}
\tau(x)x^{-1} \otimes x=1 \otimes \tau(x), \qquad x\in G_0(R). 
\end{align*}
Hence $\tau(x)$ is a scalar multiple of $x$, and we get a character  $\alpha\in\widehat{G_0(R)}$. Similarly for $\omega$; we call the corresponding character $\beta$ and \eqref{eqn:chars1} follows.

Now, we see from \eqref{eqn:id2} that there is a character $\chi\in\widehat{G_0(R)}$ such that $g\cdot x=\chi(x)x$; which shows that $g$ is a quasi-central group-like element. Hence $\tau(x)=\gamma(x)g\cdot x$, for $\gamma(x)=\alpha(x)\chi(x)^{-1}$. Analogously, $h\in G(R)$ is quasi-central and  $\omega(x)=\delta(x)h\cdot x$, for some $\delta\in\widehat{G_0}$; \eqref{eqn:chars2} follows.

Finally, we see that 
\begin{align}\label{eqn:gh1}
gh=\alpha(h)\gamma(h)^{-1}hg,
\end{align}
whence $gh(hg)^{-1}\in k^\times$ and \eqref{eqn:gh} follows.

\smallskip

$\diamondsuit${\it The bracket rule.}
Recall that, for any $t\in G_0(R)$, $Xt=\alpha(t)\,tX$ and $Yt=\beta(t)\,tY$. Also, 
observe that both $XY$ and $YX$ are $\sigma$-central, $\sigma=\tau\omega$, in $A_R$ and hence $XY-\xi\,YX\in Z_\sigma(R)$.
Now,
\begin{align*}
\mu(XY)=& XY \otimes 1 \otimes 1 - Xh\otimes h^{-1}Y\otimes 1 + Xh\otimes h^{-1}\otimes Y\\
&\quad - gY \otimes g^{-1}X \otimes 1 + gh\otimes h^{-1}Yg^{-1}X\otimes 1 - gh\otimes h^{-1}g^{-1}X\otimes Y\\
&\quad +gY \otimes g^{-1} \otimes X - gh\otimes h^{-1}Yg^{-1}\otimes X + gh\otimes h^{-1}g^{-1}\otimes XY\\
=& XY \otimes 1 \otimes 1 - \alpha(h)\,hX\otimes h^{-1}Y\otimes 1 + \alpha(h)\,hX\otimes h^{-1}\otimes Y\\
&\quad - gY \otimes g^{-1}X \otimes 1 + \beta(g^{-1})\,gh\otimes (gh)^{-1}YX\otimes 1 - gh\otimes (gh)^{-1}X\otimes Y\\
&\quad +gY \otimes g^{-1} \otimes X - \beta(g^{-1})\,gh\otimes (gh)^{-1}Y\otimes X + gh\otimes (gh)^{-1}\otimes XY.
\end{align*}
On the other hand, using \eqref{eqn:gh},
\begin{align*}
\mu(YX)=&YX\otimes 1 \otimes 1 - \beta(g)gY \otimes g^{-1}X\otimes 1 + \beta(g)gY\otimes g^{-1}\otimes X\\
&\quad -hX\otimes h^{-1}Y \otimes 1 + \alpha(h^{-1})gh\otimes (gh)^{-1}XY\otimes 1 - gh\otimes (gh)^{-1}Y\otimes X\\
&\quad +hX\otimes h^{-1} \otimes Y - \alpha(h^{-1})gh\otimes (gh)^{-1}X\otimes Y + gh\otimes (gh)^{-1}\otimes YX.
\end{align*}
As we request $XY-\xi\,YX\in R$, we see that
\begin{align}\label{eqn:chi=xi}
\alpha(h)=\xi,  \ \beta(g)=\xi^{-1},
\end{align}
and there is a $\sigma$-central skew-primitive element $c\in Z_{\sigma}(R)\cap \cP_{gh,1}(R)$, so that $XY-\xi\,YX=c$.
As a result, %$\tau(g)=\xi\, g$, 
$\tau(h)=\xi\, h$, $\omega(g)=\xi\inv\, g$. %$\omega(h)=\xi\inv\, h$.

\smallskip
$\diamondsuit${\it On the converse.}
The previous computations show that given a Hopf-Galois algebra $(R,\mu)$, if there is a quintuple $(\tau,\omega,g,h,c,\xi)$ satisfying conditions (1) to (4), then \eqref{eqn:muY-thm} defines a Hopf-Galois structure on $A$ that extends $\mu$.

\smallskip

This ends the proof of the theorem. %\hfill\qed
\epf

We write down the case $\omega=\tau^{-1}$ for completeness. We start with a remark that will be of importance further on.

\begin{rem}\label{rem:gh}
We point out that even in the case $\omega=\tau^{-1}$ we do not have $gh=hg$ as in the Hopf algebra case, cf.~Theorem \ref{thm:hopf-algebra} and Corollary \ref{cor:hopf-algebra}. Instead, we have that $g,h$ commute up to a scalar, see \eqref{eqn:gh}. In particular, 
\begin{align}\label{eqn:gh-r}
gh\cdot r=\alpha(h)\gamma(h)^{-1}\,r, \qquad r\in R.
\end{align}
\end{rem}

\begin{corollary}
Let $R=(R,\mu)$ be a Hopf-Galois algebra. Assume given a quintuple $(\tau,g,h,c,\xi)$ where
\begin{enumerate}
	\item $\tau\in\Aut(R)$ is an algebra automorphism satisfying, for any $r\in R$,
	\begin{align}
	\label{eqn:id1-cor} \tau(r)_{(1)}\ot \tau(r)_{(2)}\ot \tau(r)_{(3)}&=\tau(r_{(1)})\ot r_{(2)} \otimes r_{(3)}=\tau(r_{(1)})\ot \tau(r_{(2)}) \otimes \tau(r_{(3)} ).
	\end{align}
	\item $g,h\in G_0(R)$ are quasi-central group-like elements such that,  for any $r\in R$,
	\begin{align}
	\label{eqn:id2-cor} g\cdot r_{(1)}\otimes g\cdot r_{(2)}\ot r_{(3)}&=\tau(r_{(1)})\ot \tau(r_{(2)}) \otimes r_{(3)},\\
	\label{eqn:id3-cor} h\cdot r_{(1)}\otimes h\cdot r_{(2)}\ot r_{(3)}&=\tau^{-1}(r_{(1)})\ot \tau^{-1}(r_{(2)}) \otimes r_{(3)}. 
	\end{align}
	\item $\xi\in k^\ast$ is a scalar such that $\tau(g)=\xi\,g$ and $\tau(h)=\xi\,h$.
	\item $c\in Z(R)\cap \cP_{gh,1}(R)$ is a central skew-primitive element, i.e.
	\[
	\mu(c)=c\ot1\ot1-gh\ot (gh)^{-1}c\ot 1+gh\ot (gh)^{-1}\ot c.
	\]
\end{enumerate}
Then there is a unique Hopf-Galois algebra structure on $A=A(R,X,Y,\tau,c,\xi)$ extending $\mu$ and such that \eqref{eqn:muY-thm} holds.

Conversely, given $\tau\in\Aut(R)$, $c\in Z(R)$ and $\xi\in k^\ast$, if \eqref{eqn:muY-thm} defines a Hopf-Galois algebra structure on $A(R,X,Y,\tau,c,\xi)$ extending $\mu$, for some $g,h\in G_0(R)$, then conditions $(1)$ to $(4)$ are satisfied.
\end{corollary}
\pf
This is Theorem \ref{thm:mainmain} for $\omega=\tau^{-1}$. In particular, notice that the four conditions \eqref{eqn:id1}--\eqref{eqn:id3} become just three, as we obtain \eqref{eqn:id1-cor} and \eqref{eqn:id2-cor}, together with the additional \eqref{eqn:id3-cor}. For the remaining identity in \eqref{eqn:id1w}, namely
\begin{align*}
\tau^{-1}(r_{(1)})\ot r_{(2)} \otimes r_{(3)}=\tau^{-1}(r_{(1)})\ot \tau^{-1}(r_{(2)}) \otimes \tau^{-1}(r_{(3)} ), 
\end{align*} we remark that it follows by applying $\tau^{-2}\ot\tau^{-1}\ot\tau^{-1}$ to the corresponding identity in \eqref{eqn:id1}, as 
 $\mu\circ \tau^{-1}=(\tau^{-1})^{\ot 3}\circ\mu$ .
\epf

\section{The Hopf algebra associated to an ambiskew Hopf-Galois algebra}\label{sec:hopf-to-ambiskew}

Let $(R,\mu)$ be a Hopf-Galois algebra. Recall that this is equivalent to a structure of Hopf-Galois object over certain Hopf algebra $\underline{H}(R)$, by Theorem \ref{thm:hr}. In this section, we study the corresponding Hopf algebra $\underline{H}(A_R)$, where $A_R=A(R,X,Y,\tau,\omega,c,\xi)$ is the generalized ambiskew polynomial algebra that extends $(R,\mu)$ as in Theorem \ref{thm:mainmain}. Recall that $\sigma=\tau\omega=\omega\tau$.

We start with a technical lemma.
%\begin{lemma}\label{lem:eqn3}
%Let $(R,\mu_R)$ and $A_R$ be as above. Then
%\begin{align}\label{eqn3}
%\mu(r)=r_{(1)}\ot \tau(r_{(2)})\ot \tau(r_{(3)}).
%\end{align}
%It also follows that 
%\begin{align}\label{eqn4}
%\mu(r)=gh\cdot r_{(1)}\ot gh\cdot r_{(2)}\ot r_{(3)}.
%\end{align}
%\end{lemma}
%\pf
%By \eqref{eqn:id1} for $\sigma^{-1}$, we have
%\begin{align*}
%\tau^{-1}(r_{(1)})\ot r_{(2)}\ot r_{(3)}=\sigma^{-1}(r_{(1)})\ot \sigma^{-1}(r_{(2)})\ot \sigma^{-1}(r_{(3)}).
%\end{align*}
%Then \eqref{eqn3} follows by applying $\sigma^{\ot3}$ to this identity. Finally \eqref{eqn4} is the combination of \eqref{eqn:id2} and \eqref{eqn:id3}.
%\epf
We shall stick to the notation $r\ot s\coloneqq\sum_ir_i\ot s_i\in \underline{H}(R)$ for a generic element in $\underline{H}(R)$, see Remark \ref{rem:sweedler}.

\begin{lemma}\label{lem:H(A)}
Let $(R,\mu_R)$ and $A_R$ be as above. Then
\begin{enumerate}
\item $\underline{H}(R)\subset \underline{H}(A_R)$.
\item $(g\ot g^{-1})(h\ot h^{-1})=(h\ot h^{-1})(g\ot g^{-1})$.
\item If $r\ot s\in \underline{H}(R)$, then 
\begin{align*}
\tau(r)\ot \tau(s)&=r\ot s, & \omega(r)\ot \omega(s)&=r\ot s.
\end{align*}
\item If $r\ot s\in \underline{H}(R)$, then
\begin{align*}
g\cdot r\ot \tau(s) &\in \underline{H}(R), & h\cdot r\ot \omega(s) &\in \underline{H}(R).
\end{align*}
In particular, $(gh\cdot r)\ot \sigma(s) \in \underline{H}(R)$. 
\end{enumerate}
\end{lemma}
\pf
(1) This is clear, as $\mu_{A|R}=\mu_R$. (2) is \eqref{eqn:gh}.

(3) Recall that $r\ot s\in \underline{H}(R)$ if and only if 
\[
rs_{(1)}\ot s_{(2)}\ot  s_{(3)}=1\ot r\ot s.
\]
If we apply $\tau^{\ot 3}$ to this equality, as $\tau(1)=1$, then we get
\[
1\ot \tau(r)\ot \tau(s)=\tau (rs_{(1)})\ot \tau(s_{(2)})\ot  \tau(s_{(3)})\stackrel{\eqref{eqn:id1}}{=}\tau(rs_{(1)})\ot s_{(2)}\ot  s_{(3)}=1\ot r\ot s.
\]
The result for $\omega$ is analogous.

(4) By \eqref{eqn:id1} and \eqref{eqn:id2}:
\begin{align*}
(g\cdot r) \tau(s)_{(1)}\ot \tau(s)_{(2)}\ot \tau(s)_{(3)}&=(g\cdot r) (g\cdot s_{(1)})\ot g\cdot s_{(2)}\ot \tau(s_{(3)})=g\cdot 1\ot g\cdot r\ot \tau(s).
\end{align*} 
Thus $g\cdot r\ot \tau(s)\in \underline{H}(R)$. The other case is similar.
\epf

We use Lemma \ref{lem:H(A)} to define algebra automorphisms $\tau',\omega'\in \Aut \underline{H}(R)$ via the formulas
\begin{align}\label{eqn:tau-omega}
\tau'(r\ot s)&=g\cdot r\ot \tau(s), & \omega'(r\ot s)&=h\cdot r\ot \omega(s), \qquad r\ot s \in \underline{H}(R).
\end{align}
Notice that $\tau'\omega'=\omega'\tau'$ as $g,h$ are quasi-central and thus  $gh\cdot r=hg\cdot r$, $r\in R$. 

\begin{remark}
If $\tau\omega=\id$, that is if $gh\cdot r\ot s=r\ot s$ for all $r\ot s \in \underline{H}(R)$, then 
$gh\in Z(G_0(R))$ and $gh=hg$.
\pf
Indeed, let $t\in G_0(R)$ and fix $r\ot s=t\ot t^{-1}\in \underline{H}(R)$. Then $(gh\cdot t)\ot t^{-1}=t\ot t^{-1}$ which implies $gh\cdot t=t$, i.e.~$gh\in Z(G_0(R))$. If $t=h$, then it follows that $gh=hg$.
\epf
\end{remark}

\begin{theorem}\label{thm:H(A)}
Let $(R,\mu)$ be a Hopf-Galois algebra, let $(\tau,\omega,g,h,c,\xi)$ be a sextuple satisfying the conditions of Theorem 	\ref{thm:mainmain} and let $A_R=A(R,X,Y, \tau,\omega, c, \xi)$ be the associated Hopf-Galois algebra.

Consider $\tau',\omega'\in \Aut \underline{H}(R)$ as in \eqref{eqn:tau-omega}, the group-like elements $g^{-1}\otimes g, h^{-1}\otimes h \in G(\underline{H}(R))$, the scalar $\lambda\in k^*$ such that $gh=\lambda hg$,  and the elements $c'\in \underline{H}(R)$ and $\xi'\in k^*$ defined by
\begin{align*}
c'&=(gh)^{-1}\ot c - (gh)^{-1}c\ot 1, & \xi'&=\xi\lambda^{-1}.
\end{align*}
Then $A(\underline{H}(R),E,F, \tau',\omega', c', \xi')$ has a Hopf algebra structure as in Theorem \ref{thm:hopf-algebra}, and we have a Hopf algebra isomorphism
\[ \underline{H}(A_R)\simeq A(\underline{H}(R),E,F,\tau',\omega',c',\xi'). \]
	\end{theorem}
	
\pf
\textsl{Step 1.}
We first have to show that the conditions in Theorem \ref{thm:hopf-algebra} are satisfied. Define $\alpha, \beta : \underline{H}(R) \to k$ by $\alpha(r\otimes s)= (g\cdot r)\tau(s)$ and   $\beta(r\otimes s)= (h\cdot r)\omega(s)$. It is clear that $\alpha$ and $\beta$ are characters on $\underline{H}(R)$. We have
\begin{align*} \alpha((r\otimes s)_{(1)})(r\otimes s)_{(2)}&=\alpha(r\otimes s_{(1)})s_{(2)}\otimes s_{(3)}
=(g\cdot r)\tau(s_{(1)})s_{(2}\otimes s_{(3)}\\
&\stackrel{\eqref{eqn:id1}}{=}(g\cdot r)\tau(s_{(1)})\tau(s_{(2})\otimes \tau(s_{(3)}) 
= g\cdot r \otimes \tau(s) =\tau'(r\otimes s)
\end{align*}
so that \eqref{eqn:tau-omega-hopf} holds for $\tau'$, and for $\omega'$. Similarly
\begin{align*} g^{-1}\otimes g\cdot (r\otimes s)_{(1)}\alpha((r\otimes s)_{(2)})&=g^{-1}\otimes g \cdot (r\otimes s_{(1)})\alpha(s_{(2)}\otimes s_{(3)})
=(g\cdot r \otimes g\cdot s_{(1)}) g\cdot s_{(2)}\tau(s_{(3)})\\
&\stackrel{\eqref{eqn:id2}}{=}(g\cdot r\otimes \tau(s_{(1)}))\tau(s_{(2})\tau(s_{(3)}) 
= g\cdot r \otimes \tau(s) =\tau'(r\otimes s)
\end{align*}
and \eqref{eqn:alpha-beta} holds for $\tau'$, and for $\omega'$ as well. 

It is clear that $\tau'$ and $\omega'$ commute since $\tau$ and $\omega$ commute and $g$ and $h$ commute up to a scalar, so the characters $\alpha$ and $\beta$ convolution commute. 

We have
\[ \alpha(h^{-1}\otimes h)=g\cdot h^{-1}\tau(h) =\lambda^{-1}h^{-1}\xi h=\xi' \]
and similarly $\beta(g^{-1}\otimes g)=\xi'^{-1}$. It is an immediate verification that $c'\in \underline{H}(R)$ and that $c'$ is $((gh)^{-1}\otimes gh,1)$-primitive, with $(gh)^{-1}\otimes gh=(g^{-1}\otimes g) (h^{-1}\otimes h)$. We have finally, for $\sigma'=\tau'\omega'$,
\begin{align*}c'(r\otimes s)&= r(gh)^{-1}\otimes cs - r(gh)^{-1}c\otimes s = (gh)^{-1}(gh)\cdot r \otimes \sigma(s)c -(gh)^{-1}(gh\cdot r) c \otimes s \\
& = (gh)^{-1}(gh)\cdot r \otimes \sigma(s)c -(gh)^{-1}c\sigma^{-1}(gh\cdot r)  \otimes s \\
& =  (gh)^{-1}(gh)\cdot r \otimes \sigma(s)c -(gh)^{-1}c(gh\cdot r)  \otimes \sigma(s) \ \text{(by Lemma \ref{lem:H(A)})}\\
& = ((gh)\cdot r \otimes \sigma(s)) c' = \sigma'(r\otimes s) c' 
\end{align*}
and hence $c'$ is $\sigma'$-central.

\textsl{Step 2.} We now check that $A_R=A(R,X,Y,\tau,\omega,c,\xi)$ has a natural $A(\underline{H}(R),E,F, \tau',\omega', c', \xi')$-comodule algebra structure. We start with preliminary computations. We have
\begin{align*}
\left(g\otimes E + X\otimes 1\otimes 1\right)&r_{(1)}\otimes r_{(2)}\otimes r_{(3)}  = 
gr_{(1)}\otimes \left(E(r_{(2)}\otimes r_{(3))})\right)+ Xr_{(1)}\otimes r_{(2)}\otimes r_{(3)}\\
& = gr_{(1)}\otimes \left(g\cdot r_{(2)}\otimes \tau(r_{(3)})E\right) + \tau(r_{(1)})X\otimes r_{(2)}\otimes r_{(3)} \\
& \stackrel{\eqref{eqn:id1}}{=} (g\cdot r_{(1)})g\otimes \left(g\cdot r_{(2)}\otimes \tau(r_{(3)})E)\right) + \tau(r)_{(1)}X\otimes \tau(r)_{(2)}\otimes \tau(r)_{(3)}\\
& \stackrel{\eqref{eqn:id2}}{=} \tau(r_{(1)})g\otimes \left(\tau(r_{(2)})\otimes \tau(r_{(3)})E\right) + \tau(r)_{(1)}X\otimes \tau(r)_{(2)}\otimes \tau(r)_{(3)}\\
&  \stackrel{\eqref{eqn:id1}}{=} \tau(r)_{(1)}g\otimes \left(\tau(r)_{(2)}\otimes \tau(r)_{(3)}E\right) + \tau(r)_{(1)}X\otimes \tau(r)_{(2)}\otimes \tau(r)_{(3)}\\
& =   \tau(r)_{(1)}\otimes \tau(r)_{(2)}\otimes \tau(r_{(3)}) \left(g\otimes E + X\otimes 1\otimes 1\right)
\end{align*}
and similarly
\[  \left(h\otimes F + Y\otimes 1\otimes 1\right)r_{(1)}\otimes r_{(2)}\otimes r_{(3)}  =  \omega(r)_{(1)}\otimes \omega(r)_{(2)}\otimes \omega(r_{(3)}) \left(h\otimes F + Y\otimes 1\otimes 1\right).\]
We also have
\begin{align*} &\left(g\otimes E + X\otimes 1\otimes 1\right)\left(h\otimes F + Y\otimes 1\otimes 1\right) -\xi\left(h\otimes F + Y\otimes 1\otimes 1\right)\left(g\otimes E + X\otimes 1\otimes 1\right) \\
& = gh \otimes (EF-\xi'FE) + (XY-\xi YX)\otimes 1 \otimes 1 + (gY-\xi Yg)\otimes E + (Xh-\xi hX)\otimes F \\
&=gh\otimes c' + c\otimes 1 \otimes 1 = gh\otimes(gh)^{-1}\ot c - gh\otimes (gh)^{-1}c\ot 1+c\ot 1 \ot 1\\
& = c_{(1)}\otimes c_{(2)}\otimes c_{(3)}
\end{align*}
These computations show that we have an algebra map 
\begin{align*}
A_R &\longrightarrow A_R \otimes A(\underline{H}(R),E,F, \tau',\omega', c', \xi'),
\end{align*}
defined by
\begin{align*}
R\ni r& \longmapsto r_{(1)}\otimes r_{(2)}\otimes r_{(3)},&	X & \longmapsto g\otimes E + X\otimes 1 \otimes 1,
	& Y & \longmapsto  h \ot F + Y\ot 1\ot 1
\end{align*}
% $R\ni r \longmapsto r_{(1)}\otimes r_{(2)}\otimes r_{(3)}$ and
%\begin{align*}
%X & \longmapsto g\otimes E + X\otimes 1 \otimes 1,
%& Y & \longmapsto  h \ot F + Y\ot 1\ot 1
%\end{align*}
which is easily seen to furnish the expected comodule algebra structure. 

\textsl{Step 3.} We now prove that the previous comodule algebra structure makes $A_R$ a Hopf-Galois object over $A(\underline{H}(R))=A(\underline{H}(R),E,F, \tau',\omega', c', \xi')$. For this, we claim that there exists an algebra map 
\begin{align*}
\theta : A(\underline{H}(R)) &\longrightarrow A_R^{\rm op} \otimes A_R  \\
r\otimes s, E, F & \longmapsto r\otimes s, E':=g^{-1}\ot X - g^{-1}X\ot 1, F':=h^{-1}\ot Y -h^{-1}Y\ot 1
\end{align*} 
such that $r\otimes s\longmapsto r\otimes s$  and 
\begin{align*}
E & \longmapsto E':=g^{-1}\ot X - g^{-1}X\ot 1, & F&\longmapsto F':=h^{-1}\ot Y -h^{-1}Y\ot 1.
\end{align*} 
To prove the claim, we have to check a number of identities. We have, for $r\ot s \in \underline{H}(R)$, 
\begin{align*}
E'(r\ot s)&=rg^{-1}\ot Xs - rg^{-1}X \ot s=rg^{-1}\ot \tau(s)  X - rg^{-1}X\ot s\\
&=(g\cdot r\ot \tau(s))
(g^{-1}\ot X) -
(\tau^{-1}(g\cdot r)\ot s) (g^{-1}X\ot 1)\\
&=(g\cdot r\ot \tau(s))(g^{-1}\ot X - g^{-1}X\ot 1) \\
& =\tau'(r\ot s) E'
\end{align*}
where the third equality holds since $g\cdot r\ot \tau(s)=\tau^{-1}(g\cdot r)\ot s$, by Lemma \ref{lem:H(A)}.
Similarly $F'(r\ot s)=\omega'(r\ot s)F'$.
Finally we see that:
\begin{align*}
E'F'&=h^{-1}g^{-1}\ot XY - h^{-1}Yg^{-1}\ot X   -    h^{-1}g^{-1}X\ot Y + h^{-1}Yg^{-1}X\ot 1\\
&=(gh)^{-1}\ot XY - \xi(gh)^{-1}Y\ot X - (gh)^{-1}X\ot Y + \xi(gh)^{-1}YX\ot 1.\\
F'E'&=\lambda^{-1}(gh)^{-1}\ot YX - \lambda^{-1}\xi^{-1} (gh)^{-1}X\ot Y   -  \lambda^{-1}(gh)^{-1}Y\ot X + \lambda^{-1}\xi^{-1} (gh)^{-1}XY\ot 1.
\end{align*}
Hence we have
\begin{align*}
E'F'-\xi'F'E'&=(gh)^{-1}\ot (XY-\xi YX)-(gh)^{-1}(XY-\xi YX)\ot 1\\
&=(gh)^{-1}\ot c-(gh)^{-1}c\ot 1=c'.
\end{align*}
We have thus proved that the announced algebra map is well-defined. It is now an immediate verification that the composite map
\[  A_R\ot A(\underline{H}(R))  \overset{{\rm id} \ot \theta} \longrightarrow  A_R \ot A_R\ot A_R \overset{{\rm mult}\ot {\rm id}} \longrightarrow A_R \ot A_R\]
is inverse to the Galois map in Section \ref{sec:hopf-galois}, so we conclude that $A_R$ is a Hopf-Galois object over $A(\underline{H}(R))$.

\textsl{Step 4.} To finish the proof, first notice that the elements $E'=g^{-1}\ot X - g^{-1}X\ot 1$ and $F'=h^{-1}\ot Y -h^{-1}Y\ot 1$ from the previous step belong to $\underline{H}(A_R)$, and that the computations done in this previous step ensure the existence of a Hopf algebra map 
\begin{align*}
f: A(\underline{H}(R)) & \longrightarrow \underline{H}(R),
\end{align*}
given by
\begin{align*}
	r\otimes s&\longmapsto r \ot s, & E&\longmapsto E', & F&\longmapsto F'.
\end{align*}
It is obvious that $f$ commutes with the respective coactions on $A_R$, so using the bijective canonical maps of our two Hopf-Galois objects, we conclude that $f$ is an isomorphism.
\epf	

While the group-like and skew-primitive elements in a Hopf-Galois algebra $R$ can be described from those of $\underline{H}(R)$ thanks to Lemma \ref{lem:group}, the proof of the previous theorem also provides the description of the automorphisms in Theorem \ref{thm:mainmain} in terms of characters of $\underline{H}(R)$:
	
\begin{corollary}\label{cor:expli-tau}
	Let $(R,\mu)$ be Hopf-Galois algebra, let $(\tau,\omega,g,h,c,\xi)$ be a sextuple satisfying the conditions of Theorem 	\ref{thm:mainmain}. Then there are convolution commuting characters $\alpha, \beta : \underline{H}(R) \to k$ such that for $r \in R$,%group-like elements $g',h' \in G(\underline{H}(R))$, a skew-primitive $c'\in \mathcal P_{g'h',1}(\underline{H}(R))$ satisfying the following conditions:
 \[\tau(r) = \alpha(r_{(2)}\otimes r_{(3)}) g\cdot r_{(1)}, \quad \omega(r) = \beta(r_{(2)}\otimes r_{(3)}) h\cdot r_{(1)}
		 \]
and for $r\otimes s \in \underline{H}(R)$,
\[ \alpha((r\otimes s)_{(1)})(r\otimes s)_{(2)}= g^{-1}\otimes g\cdot (r\otimes s)_{(1)}\alpha((r\otimes s)_{(2)}),\]
\[    \beta((r\otimes s)_{(1)})(r\otimes s)_{(2)}= h^{-1}\otimes h\cdot (r\otimes s)_{(1)}\beta((r\otimes s)_{(2)}).\]
\end{corollary}

\pf
We have, for any $r\in R$, by the construction of $\tau'$,
\[ r_{(1)}\ot \tau'(r_{(2)}\ot r_{(3)}) = r_{(1)} \ot g\cdot r_{(2)}\ot \tau(r_{(3)})\]
hence
\[ g\cdot r_{(1)}\ot \tau'(r_{(2)}\ot r_{(3)}) = g\cdot r_{(1)} \ot g\cdot r_{(2)}\ot \tau(r_{(3)}).\]
Multiplying everything, we get, by the definition of $\alpha$ in the previous proof the announced identity for $\tau$ and $\alpha$, and similarly for $\omega$ and $\beta$. The final identities have been proved in the previous proof as well.
\epf

Notice that when $\omega=\tau\inv$, that is when $A_R=A(R,X,Y,\tau,c,\xi)$ is an ambiskew polynomial algebra, for some $\tau\in\Aut R$, then we do not necessarily have $\omega'=\tau^{'-1}$, as we get 
\[
\tau'\omega'(g^{-1}\ot g)=gh\cdot g^{-1}\ot g=\lambda g^{-1}\ot g.
\]
%see \eqref{eqn:gh-r}: here $\alpha,\gamma\in\widehat{G(R)}$ are such that $\sigma(x)=\alpha(x)x=\gamma(x)\,g\cdot x$, $x\in G(R)$. 
Thus it becomes evident that we need to work with generalized ambiskew polynomial algebras rather than plain ambiskew polynomial algebras. We resume this situation in the next corollary, see \S\ref{sec:illustration} for a concrete example.
\begin{corollary}\label{coro:notid}
Let $(R,\mu)$ be a Hopf-Galois algebra and assume that $A_R=A(R,X,Y,\tau,c,\xi)$ is a Hopf-Galois algebra. Then the corresponding Hopf algebra $\underline{H}(A_R)$ is a generalized ambiskew polynomial algebra $A(\underline{H}(R),\tau',\omega',c',\xi')$, where $\tau',\omega'\in\Aut \underline{H}(R)$ are commuting automorphisms not necessarily inverse to each other.
\end{corollary}

\section{Examples}\label{sec:examples}

\subsection{Baby example: $R=k$} It is immediate that Theorem \ref{thm:mainmain} produces only two examples: the Hopf algebra $k[X,Y]$ with primitive $X,Y$, and the Weyl algebra $A_1(k)$ as in Example \ref{ex:Weyl}.

\subsection{Almost general example: cleft Galois objects} \label{sub:cleft} Before going into other specific examples, we examine the case of cleft Galois objects, i.e. those obtained by deforming the multiplication of a Hopf algebra by a Hopf $2$-cocycle. This is the most studied type of Hopf-Galois objects, but there exist some that are not of this type, see \cite{bic}.

Let $H$ be a Hopf algebra and let $\pi : H \otimes H \to k$ be a $2$-cocycle on $H$, i.e. $\pi$ is convolution invertible, $\pi(x,1)=\varepsilon(x)=\pi(1,x)$ for any $x\in H$, and we have, for any $x,y,z\in H$,
\[ \pi(x_{(1)}, y_{(1)})\pi(x_{(2)}y_{(2)},z)= \pi(y_{(1)},z_{(1)})\pi(x,y_{(2)}z_{(2)}).\]
The algebra $_{\pi}\!H$ has $H$ as underlying vector space, and product defined by
\[ x . y = \pi(x_{(1)},y_{(1)}) x_{(2)}y_{(2)}.\]
It follows from standard cocycle identities \cite{doi} that $_{\pi}\!H$ is a Hopf-Galois algebra with
\begin{align*}
\mu :{ _{\pi}\!H} &\longrightarrow {_{\pi}\!H}\ot {_{\pi}\!H^{\rm op}} \ot {_{\pi}\!H}  \\
x &\longmapsto \pi^{-1}(S(x_{(3)}), x_{(4)})x_{(1)} \ot S(x_{(2)}) \ot x_{(5)}
\end{align*}
and that the map
\begin{align*}
H &\longrightarrow \underline{H}({ _{\pi}\!H}) \subset {_{\pi}\!H^{\rm op}} \ot {_{\pi}\!H}\\
x & \longmapsto  \pi^{-1}(S(x_{(2)}), x_{(3)})S(x_{(1)}) \ot x_{(4)})
\end{align*}
is a Hopf algebra isomorphism. We thus see from Lemma \ref{lem:group} that 
\[ G_0({ _{\pi}\!H}) = \{\lambda g, \ \lambda\in k^*, g\in G(H)\}\]
and for $g,h \in G({H})$, 
\[ \mathcal P_{g,h}^0({ _{\pi}\!H})=\{ \lambda1 + x, \lambda \in k, x\in \mathcal P_{g,h}(H)\}. \]
Using this description and Corollary \ref{cor:expli-tau}, we then see that the sextuples $(\tau,\omega,g,h,c,\xi)$ as in Theorem \ref{thm:mainmain} can be chosen as follows:
\begin{enumerate}
	\item $g,h$ are central group-like elements in $H$; 
	\item  the automorphisms $\tau, \omega  \in {\rm Aut} { _{\pi}\!H}$  are of the form
	\[ \tau(x)= \alpha(x_{(2)}) g\cdot x_{(1)}, \ \omega(x)= \beta(x_{(2)})h\cdot x_{(1)}\]
	for some convolution commuting algebra maps $\alpha, \beta : H \to k$ that also satisfy 
	\[  \alpha(x_{(2)}) x_{(1)} = \alpha({x_{(1)}}) g\cdot x_{(2)} , \   \beta(x_{(2)})x_{(1)}= \beta(x_{(1)})h\cdot x_{(2)}\]
	(these last identities are in $H$) and  $\alpha(h)=\beta(g)^{-1}$, $\xi= \alpha(h) \lambda$, where $\lambda$ is such that $gh=\lambda hg$ in ${ _{\pi}\!H}$.
	\item  $c$ is $\sigma$-central for $\sigma=\tau\omega$, and $c=\gamma 1+x$, for $\gamma\in k$ and $x \in  \mathcal P_{gh,1}(H)$;
\end{enumerate}

\subsection{The example $R=k[T]$}
We now examine the case when $R= \underline{R}(k[T])$. The analysis of the previous subsection applies: we have $G_0(R)=k^*$ and $\cP_{1,1}^0(R) =k+kT$ and  for the sextuples $(\tau,\omega,g,h,c,\xi)$ of Theorem \ref{thm:mainmain}, the automorphisms $\tau$, $\omega$  are of the form
$\tau(T)=T + \alpha 1$, $\omega(T)=T +\beta 1$ for some $\alpha, \beta\in k$. 
Moreover we must have $\beta =-\alpha$ if $c\not=0$ because of the $\sigma$-centrality of $c$.

At the end, after suitable changes in the generators, we find that a Hopf-Galois algebra $A= A(R,\tau,\omega,g,h,c,\xi)$ as in  Theorem \ref{thm:mainmain} belongs to one of the following classes.
\begin{enumerate}
	\item  $A$ is the enveloping Hopf algebra of the Lie algebra $\mathfrak{sl}_2$:
	\[  A=U(\mathfrak{sl}_2) =k\langle T, X, Y \ | \ [X,T]=X, [Y,T]=-Y, \ [X,Y]=T \rangle.  \]
\item $A$ is the enveloping Hopf algebra of the $3$-dimensional Heisenberg Lie algebra $\mathfrak{h}$:
\[  A=U(\mathfrak{h}) =k\langle T, X, Y \ | \ [X,T]=0= [Y,T], \ [X,Y]=T \rangle.  \]
\item $A=U_\beta$, $\beta\in k$, is the enveloping algebra of a solvable $3$-dimensional Lie algebra:
\[  A=U_\beta =k\langle T, X, Y \ | \ [X,T]=X, \ [Y,T] =\beta Y, \ [X,Y]=0 \rangle.  \]
\item $A=k[T,X,Y]$ is abelian.
\item $A$ is a Hopf-Galois algebra and not a Hopf algebra, with 
\[  A =U_{\alpha, \lambda}= k\langle T, X, Y \ | \ [X,T]=\alpha X, [Y,T]=-\alpha Y, \ [X,Y]=1+\lambda T \rangle, \ \alpha\in \{0,1\}, \ \lambda\in k.  \]
Moreover, the Hopf algebra $\underline{H}(U_{\alpha,\lambda})$ is 
\begin{enumerate}
	\item $U(\mathfrak{sl}_2) $ if $\alpha=1$ and $\lambda \not=0$;
	\item  $U_{-1}$ if $\alpha=1$ and $\lambda =0$;
	\item  $U(\mathfrak{h})$ if $\alpha=0$ and $\lambda \not=0$;
	\item  $k[T,X,Y]$ if $\alpha=0$ and $\lambda =0$.
\end{enumerate}
\end{enumerate}

\subsection{The example $R=k[T,T^{-1}]=k\mathbb Z$}
In this part, we set $R= \underline{R}(H)$, for $H=k\mathbb Z=k\lg K^{\pm1}\rg$. As above, we choose two (central) group-like elements $g=K^n$, $h=K^m$, $m,n\in\mathbb Z$. We fix two characters $\alpha,\beta\colon H\to k$, which are determined by nonzero scalars $\alpha\coloneqq\alpha(K)$, $\beta\coloneqq\beta(K)$, satisfying $\xi=\alpha^m=\beta^{-n}$ (as $\lambda=1$). If $c\neq 0$, then there is $\gamma,\mu\in k$ such that $c=\gamma+\mu(1-gh)$. As $c$ is $\sigma$-central, it follows that in this case $\beta=\alpha^{-1}$ (hence $\alpha^n=\alpha^m=\xi$). 

Thus, a Hopf-Galois algebra $A= A(R,\tau,\omega,g,h,c,\xi)$ as in Theorem \ref{thm:mainmain}, with $g=K^n$ and $h=K^m$, belongs to one of the following classes.
\begin{enumerate}
	\item  For each $\alpha,\beta$ with $\xi=\alpha^m=\beta^{-n}$:
	\[  A=k\langle K, X, Y \ | \ XK=\alpha KX, YK=\beta KY, \ XY-\xi\,YX=0 \rangle. \]
	\item   For each $\alpha$ with $\xi=\alpha^m=\alpha^{n}$ and $\gamma,\mu\in k$:
	\[    A=k\langle K, X, Y \ | \ XK=\alpha KX, YK=\alpha^{-1} KY, \ XY-\xi\,YX=\gamma+\mu(1-K^{n+m}) \rangle.\]
\end{enumerate}
In the last class, we recognize the Hopf algebra $U_q(\sl_2)$ when $n=m=-1$, $\gamma=0$, $\mu=\frac{1}{q-q^{-1}}$ and $\xi=q^2=\alpha^{-1}$. 
On the other hand, when $\gamma\neq 0$, we obtain the $U_q(\sl_2)$-cleft objects $A_{(a)}$ as computed in \cite[Lemma 16]{g}, for $a=\gamma+\mu$ and $n=m=-1, \mu=\frac{1}{q-q^{-1}}$ as before.

\subsection{The example $R=k\langle T_{1}^{\pm 1}, T_2^{\pm 1} | T_1T_2 = q T_2T_1\rangle$}\label{sec:illustration} Here $q \in k^*$ and $R={}_\pi H$ is a twisted group algebra of the group algebra $H=k\mathbb Z^2=k\lg K_1^{\pm 1}, K_2^{\pm 1} | K_1K_2=K_2K_1 \rg$ and $\pi$ such that $\frac{\pi(K_1,K_2)}{\pi(K_2,K_1)}=q$. Subsection \ref{sub:cleft} applies to give the complete description of all the generalized ambiskew polynomial algebras that are Hopf-Galois as in Theorem \ref{thm:mainmain}. We focus on giving an example to illustrate Corollary \ref{coro:notid}, see \S\ref{sec:illustration-case1} below. In particular, we assume that $q\neq 1$. 

We choose $g=K_1, h=K_2\in G(H)$, so $gh=q\,hg$ in $R$, that is, $\lambda=q$. Now, we get that
\begin{align*}
\tau(T_1)&= \alpha(K_1)T_1, & \tau(T_2)&= q\alpha(K_2)T_2, & \omega(T_1)&= q^{-1}\beta(K_1)T_1, & \omega(T_2)&= \beta(K_2)T_2,
\end{align*}
for some algebra maps $\alpha, \beta : H \to k$ such that $\alpha(K_2)=\xi q^{-1}=\beta(K_1)^{-1}$.

In particular, $\sigma(T_1)=\xi^{-1} \alpha(K_1)T_1$ and $\sigma(T_2)=\xi\beta(K_2)T_2$.

Now $c=\gamma+\mu(1-T_1T_2)$. If $c\neq 0$, then being $\sigma$-central amounts to
\begin{enumerate}
\item If $\mu=0$, then $\xi=\alpha(K_1)=\beta(K_2)^{-1}$. Thus, $\omega=\tau^{-1}$ and hence $\sigma=\id$.
\item If $\mu\neq 0$ and $\gamma=-\mu$, then $\alpha(K_1)=q^{-1}\xi$ and $\beta(K_2)=q\xi^{-1}$. Hence $\sigma\not=\id$.
%\sigma(K_1)=q\,K_1$, $\sigma(K_2)=q\,K_2$.
\item Otherwise, $q=1$ and this is a contradiction.
\end{enumerate}

Next we analyze cases (1) and (2).
\subsubsection{Case \emph{(1)}}\label{sec:illustration-case1}
Fix $\mu=0$, $\gamma=1$, so $\tau\omega=\id$ and $c=1$. 
This will provide an example for Corollary \ref{coro:notid}.
We have
\begin{align*}
\tau(T_1)&=\xi T_1, & \tau(T_2)&=\xi T_2, & \omega(T_1)&=\xi^{-1} T_1, & \omega(T_2)&=\xi^{-1} T_2.  
\end{align*}
Let
$A_R=A(R,X,Y,\tau,c,\xi)$ be the corresponding Hopf-Galois algebra. This is generated by $T_1^{\pm1}$, $T_2^{\pm1}$, $X$ and $Y$ with relations
\begin{align*}
T_1T_2&=qT_2T_1, & XT_1&=\xi T_1X, & XT_2&=\xi T_2X, & YT_1&=\xi^{-1} T_1Y, & YT_2&=\xi^{-1} T_2Y, & XY-\xi YX=1.
\end{align*}

Now, by Theorem \ref{thm:H(A)}, the Hopf algebra $\underline{H}(A_R)$ is a generalized ambiskew polynomial algebra $A(\underline{H}(R),\tau',\omega',c',\xi')$. 

In this case, $\underline{H}(R)\simeq H$ via the assignment $T_1^{-1}\ot T_1\rightsquigarrow K_1$,  $T_2^{-1}\ot T_2\rightsquigarrow K_2$. 
Recall that 
\begin{align*}
\xi'&=\xi\lambda^{-1}=q^{-1}\xi, & c'&=(gh)^{-1}\ot c - (gh)^{-1}c\ot 1=0.
\end{align*}
Finally, $\tau',\omega'\in\Aut \underline{H}(R)$ are the automorphisms given by
\begin{align*}
\tau'(T_1^{-1}\ot T_1)&=T_1\cdot T_1^{-1}\ot \tau(T_1)=\xi T_1^{-1}\ot T_1, \\ 
\omega'(T_1^{-1}\ot T_1)&=T_2\cdot T_1^{-1}\ot \tau^{-1}(T_1)=q\xi^{-1} T_1^{-1}\ot T_1,\\
\tau'(T_2^{-1}\ot T_2)&=T_1\cdot T_2^{-1}\ot \tau(T_2)=\xi q^{-1} T_2^{-1}\ot T_2, \\
 \omega'(T_2^{-1}\ot T_2)&=T_2\cdot T_2^{-1}\ot \tau^{-1}(T_2)=\xi^{-1}T_2^{-1}\ot T_2.
\end{align*}
In other words, under  the identification above, we have that
\begin{align*}
\tau'(K_1)&=\xi\, K_1, & \tau'(K_2)&=\xi q^{-1}\, K_1, & \omega'(K_1)&=\xi^{-1}q\, K_1, & \omega'(K_2)&=\xi^{-1}\, K_2. 
\end{align*}
In particular, we see that $\tau'$ and $\omega'$ are not inverse to each other and the Hopf algebra  $\underline{H}(A_R)$ is generated by 
$K_1^{\pm1}$, $K_2^{\pm1}$, $E$ and $F$ with relations
\begin{align*}
K_1K_2&=K_2K_1, & EK_1&=\xi K_1E, & EK_2&=\xi q^{-1} K_2E,  & FK_1&=\xi^{-1}q K_1F, & FK_2&=\xi^{-1} K_2F,
\end{align*}
\begin{align*}
EF-q^{-1}\xi FE=0.
\end{align*}

\subsubsection{Case \emph{(2)}}
We have $c=\gamma T_1T_2$ and 
\begin{align*}
\tau(T_1)&=q^{-1}\xi T_1, & \tau(T_2)&=\xi T_2, & \omega(T_1)&=\xi^{-1} T_1, & \omega(T_2)&=q\xi^{-1} T_2.  
\end{align*}
Hence $A_R=A(R,X,Y,\tau,c,\xi)$ is generated by $T_1^{\pm1}$, $T_2^{\pm1}$, $X$ and $Y$ with relations
\begin{align*}
T_1T_2&=q T_2T_1, & XT_1&=q^{-1}\xi T_1X, & XT_2&=\xi  T_2X, & YT_1&=\xi^{-1} T_1Y, & YT_2&=q\xi^{-1} T_2Y,
\end{align*}
\begin{align*}
XY-\xi YX=\gamma T_1T_2.
\end{align*}
For the Hopf algebra $\underline{H}(A_R)$, we get $c'=\gamma(1-K_1K_2)$ and 
\begin{align*}
\tau'(K_1)&=q^{-1}\xi\, K_1, & \tau'(K_2)&=\xi q^{-1}\, K_1, & \omega'(K_1)&=\xi^{-1}q\, K_1, & \omega'(K_2)&=q\xi^{-1}\, K_2. 
\end{align*}
In this case, $\omega'\tau'=\id$. The Hopf algebra is now generated by 
$K_1^{\pm1}$, $K_2^{\pm1}$, $E$ and $F$ with relations
\begin{align*}
K_1K_2&=K_2K_1, & EK_1&=q^{-1}\xi K_1E, & EK_2&=\xi q^{-1} K_2E, & FK_1&=\xi^{-1}q K_1F, & FK_2&=q\xi^{-1} K_2F,
\end{align*}
\begin{align*}
EF-q^{-1}\xi FE=\gamma(1-K_1K_2).
\end{align*}

\subsection{The example $R=A_1(k)$}
Here $R=k\lg U,V | UV-VU=1\rg$ is a cleft object for the polvnomial Hopf algebra $H=k[u,v]$, with $G(H)=\{1\}$ and $\cP_{1,1}(H)=k u +k v$. As above, we obtain that the automorphisms $\tau,\omega$ are determined by algebra maps $\alpha,\beta\colon H\to k$, as $\tau(U)=U+\alpha(u)$, $\omega(U)=U+\beta(u)$, similarly for $\tau(V),\omega(V)$. We have that $c=\gamma+c_1U+c_2V$, for some $\gamma,c_1,c_2\in k$. The fact that $c$ is $\sigma=\tau\omega$-central determines that $c_1=\alpha(v)+\beta(v)$ and $c_2=-\alpha(u)-\beta(u)$.  Therefore $A_R= A(R,\tau,\omega,g,h,c,\xi)$  with $g=h=1$ is determined by scalars
$\gamma,\alpha_u,\alpha_v,\beta_u,\beta_v\in k$ and relations
\[
[U,V]=1, \qquad
[X,U]=\alpha_uX, \qquad 
[X,V]=\alpha_vX, \qquad 
[Y,U]=\beta_uY, \qquad 
[Y,V]=\beta_vY
\]
\begin{align*}
[X,Y]=\gamma+(\alpha_v+\beta_v)U-(\alpha_u+\beta_u)V.
\end{align*}If $c=\gamma+(\alpha_v+\beta_v)U-(\alpha_u+\beta_u)V=0$, then this is a Hopf algebra. When $c\neq 0$, if we identify $\underline{H}(R)$ with $H=k[u,v]$ via $1\ot U-U\ot 1\rightsquigarrow u$, $1\ot V-V\ot 1\rightsquigarrow v$ then for the corresponding Hopf algebra $\underline{H}(A_R)$ it follows that $\xi'=\xi$, $\tau'=\tau$ and $\omega'=\omega$, in such a way that $\underline{H}(A_R)=A(H,\tau,\omega,c',\xi)$ for $c'=(\alpha_v+\beta_v)u-(\alpha_u+\beta_u)v$.

\bibliographystyle{amsalpha}

\end{document}